\documentclass[a4paper,12pt]{article}
\usepackage{graphicx}
\usepackage[colorlinks=true,linkcolor=blue,citecolor=blue,urlcolor=blue]{hyperref}
\usepackage{algorithm}
\usepackage{algpseudocode}
\usepackage{amsmath,amssymb,amsfonts,bm,url}
\usepackage{subcaption}
\usepackage{booktabs}
\usepackage{multirow}
\usepackage{cite}

\title{A Unified Column Generation and Elimination Method for Solving Large-Scale Set Partitioning Problems}
\author{Yasuyuki Ihara\\
Innovation Laboratories\\
NEC Solution Innovators, Ltd.\\
1-18-7 Shinkiba, Koto-ku, Tokyo, 136-8627 Japan\\
\texttt{ihara\_mxk[ ]nec.com}\thanks{Replace [ ] with @.}}

\date{}

\begin{document}
\maketitle

\begin{abstract}
The Set Partitioning Problem is a combinatorial optimization problem with wide-ranging applicability, used to model various real-world tasks such as facility location and crew scheduling.
However, real-world applications often require solving large-scale instances that involve hundreds of thousands of variables. Although the conventional Column Generation method is popular for its computational efficiency, it lacks a guarantee for exact solutions.
This paper proposes a novel solution method integrating relaxation of Column Generation conditions and automatic elimination of redundant columns, aimed at overcoming the limitations of conventional Column Generation methods in guaranteeing exact optimal solutions.
Numerical experiments using actual bus route data reveal that while the traditional method achieves an exact solution rate of only about 3\%, the proposed method attains a rate of approximately 99\% and remarkably improves solution accuracy.
\end{abstract}

\section{Introduction}

The Set Partitioning Problem is a combinatorial optimization problem of selecting subsets from a given collection to form a partition (complete and disjoint coverage) of the original set at minimum total 
cost~\cite{SPP1,SPP2_SCP}. This problem finds diverse practical applications, including crew scheduling~\cite{CrewScheduling1,CrewScheduling2}, delivery planning~\cite{DeliveryPlanning1,DeliveryPlanning2_ColumnGeneration1}, facility location~\cite{FacilityLocation1}, bin packing~\cite{BinPacking1}, and modularity (density) maximization~\cite{ModularityMaximization1,ModularityDensityMaximization1_ColumnGeneration2}.

In many real-world applications, cases arise in which the number of variables reaches hundreds of thousands. For these large-scale instances of the Set Partitioning Problem, the Column Generation method~\cite{DeliveryPlanning2_ColumnGeneration1,ModularityDensityMaximization1_ColumnGeneration2,ColumnGeneration1,ColumnGeneration2,ColumnGeneration3,Application_ColumnGeneration1,Survey_ColumnGeneration1} has been widely adopted for its computational efficiency. Originally proposed as a solution method for large-scale linear programming problems~\cite{ColumnGeneration1,ColumnGeneration2}, Column Generation has recently been extended to address large-scale combinatorial optimization problems. In this approach, the dual of a Restricted Master Problem (RMP)—which contains only a subset of the variables from the linear relaxation—is solved, and the resulting dual solution is used to iteratively generate new promising columns based on reduced cost information. However, although this approach offers advantages in computational efficiency and ease of implementation, they typically fail to guarantee exact optimal solutions and often yield only approximate solutions. Compared to the related Set Covering Problem~\cite{SPP2_SCP}, the Set Partitioning Problem presents stricter equality constraints rather than inequality constraints, further complicating the attainment of exact solutions. Alternatively, the integration of Column Generation with branch-and-bound, known as branch-and-price~\cite{BranchPrice1,BranchPrice2,BranchPrice3}, can secure an exact solution; however, it entails complex branching procedures and increased memory usage.

To overcome these limitations, this paper introduces a new approach, termed the "Column Generation and Elimination Method," which extends traditional Column Generation techniques. By relaxing conditions within Column Generation, the proposed method significantly enhances the rate of obtaining exact optimal solutions and improves the accuracy of approximate solutions. Additionally, incorporating an automatic elimination mechanism for redundant columns mitigates the growth in problem size. Unlike branch-and-price methods, the proposed approach does not require complex branching operations and remains similarly straightforward to implement as standard Column Generation.

Numerical experiments conducted using real-world bus route data available in GTFS format demonstrate the effectiveness of the proposed approach. For the minimum crew scheduling problem—a specific case of Set Partitioning Problem with uniform cost coefficients—the exact solution rate improved dramatically from approximately 3\% using conventional Column Generation method to approximately 99\% with the proposed method. Furthermore, although both conventional and proposed methods faced difficulties obtaining exact solutions for scenarios with non-uniform cost coefficients, the proposed method achieved a notable improvement in the accuracy of approximate solutions compared to conventional methods.

The remainder of this paper is organized as follows. Section~\ref{sec:RelatedWorks} provides an overview of the Set Partitioning Problem and the conventional Column Generation method. Section~\ref{sec:ProposedMethod} describes the details of the proposed Column Generation and Elimination Method, while Section~\ref{sec:Evaluation} presents the results of numerical experiments. Finally, Section~\ref{sec:Conc} concludes the paper.


\section{Related Works}\label{sec:RelatedWorks}

\subsection{Set Partitioning Problem}\label{subsec:setpartition}

This subsection defines the Set Partitioning Problem and presents its formulation as a combinatorial optimization problem. Consider a set $ V = \{ v_i \}_{i \in I} $ along with a family of its subsets $\{ U_j \}_{j \in J}$, where $I$ and $J$ are index sets. For each $ j \in J $, let $ c_j > 0 $ denote the cost associated with selecting subset $ U_j $. 
The Set Partitioning Problem involves finding a subset $ K \subseteq J $ such that the corresponding family of subsets $ \{ U_{j_k} \}_{j_k \in K} $ forms a partition of the set $ V $, while minimizing the total cost $ \sum_{j_k \in K} c_{j_k} $. Here, a partition satisfies the following conditions:
\begin{itemize}
\item $ V = \bigcup_{j \in K} U_j $ \quad (complete coverage),
\item For any distinct $ j_1, j_2 \in K $, $ U_{j_1} \cap U_{j_2} = \emptyset $ \quad (disjoint).
\end{itemize}
Here, a binary variable $x_j$ is introduced to represent whether subset $U_j$ is included in this partition:
\begin{equation}
x_j =
\begin{cases}
1, & \text{if } j \in K, \\
0, & \text{otherwise.}
\end{cases}
\end{equation}
Then, the problem is formulated as the following combinatorial optimization 
problem~\cite{SPP1,SPP2_SCP}:
\begin{equation}
\begin{split}
& \underset{{\bm x}}{\operatorname{min}} \; {\bm c}^\top {\bm x}, \\
& \text{subject to} \quad A {\bm x} = {\bm 1}_M.
\end{split}
\end{equation}
In this formulation, the symbols are defined as follows:
\begin{itemize}
\item $ M $: the number of elements in $ V $; $ N $: the number of indices in $ J $,
\item $ {\bm x} = (x_1, \ldots, x_N)^\top $ and $ {\bm c} = (c_1, \ldots, c_N)^\top $,
\item $ {\bm 1}_M $: an $ M $-dimensional vector with all components equal to 1,
\item $ A = (a_{i,j})_{1 \leq i \leq M, \, 1 \leq j \leq N} $ is an $ M \times N $ matrix with entries defined by
\begin{equation}
a_{i,j} =
\begin{cases}
1, & \text{if } v_i \in U_j, \\
0, & \text{otherwise.}
\end{cases}
\end{equation}
\end{itemize}

\subsection{Column Generation Method}

In combinatorial optimization problems (0-1 integer programming) aimed at minimizing an objective function, the Column Generation method first solves the RMP for the linear relaxation using only a subset of the 
original column set $J$. Next, the dual solution $\bm{y}=(y_1,\ldots,y_M)^\top$ is employed to generate
and add the column which reduced cost,
\begin{equation}
\widehat{c}_j = c_j - \sum_{i=1}^M a_{ij} y_i,
\end{equation}
is negative. This process is repeated until no columns with a negative reduced cost can be found. The introduction of such columns reinforces the constraints in the dual of the RMP, thereby lowering the RMP's optimal value.

The initial column set $J_{\rm ini}$ is determined based on a feasible solution obtained using a greedy algorithm~\cite{Greedy}. Details of the greedy algorithm are provided in \ref{sec:greedy}. Although the greedy algorithm does not always yield a feasible solution, a feasible solution is obtained in most cases when the number of columns $N$ is sufficiently larger than the number of rows $M$.

The pseudocode for applying the Column Generation method to the Set Partitioning Problem until a final solution is reached is presented in {\bf Algorithm}~\ref{alg:ColGen}. In this pseudocode, the vectors $\bm{c}$ and $\bm{x}$ and the columns of the matrix $A$ restricted to the subset $K$ are denoted by $\bm{c}_K$, $\bm{x}_K$, and $A[:,K]$, respectively.

\begin{algorithm}
\caption{Column Generation Method}
\label{alg:ColGen}
\begin{algorithmic}[1]
\State \textbf{Step 1:} Employ a greedy algorithm to compute a feasible solution $\bm{z}=(z_1,\ldots,z_N)^\top$ for the original Set Partitioning Problem. From the index set $J=\{1,\ldots,N\}$, define
$K_{\rm ini}:=\{\,j\in J \mid z_j=1\,\}$
as the initial column set, and set $K \gets K_{\rm ini}$. The set $K$ will be used for column selection in the master problem.

\While{True}
   \State \textbf{Step 2a:} Solve the dual problem of the linear relaxation of the RMP restricted to the index set $K$. Specifically, solve
    \begin{equation}
    \begin{split}
    \widehat{\bm{y}} =\ & \underset{\bm{y}=(y_1,\ldots,y_M)^\top}{\operatorname{argmax}}\; \sum_{i=1}^{M} y_i,\\[1ex]
    & \text{subject to}\quad A[:,K]^\top\, \bm{y}\le \bm{c}_K,
    \end{split}
    \label{DRMP1}
    \end{equation}
    where each $y_i$ is an unconstrained continuous variable and $\bm{c}_K=(c_j)_{j\in K}$.
    
    \State \textbf{Step 2b:} Using the dual solution $\widehat{\bm{y}}$, compute 
     the reduced cost
    \begin{equation}
    \widehat{c}_j:= c_j-\widehat{\bm{y}}^\top A[:,j]\quad(j\in J\setminus K)
    \end{equation}
    for each $j$, and determine $\alpha:=\min\{\widehat{c}_j\mid j\in J\setminus K\}$.

    \State \textbf{Step 2c (Iteration Termination Check):} Define the candidate column set
    \begin{equation}
    K_{\rm gen}:=\{\,j\in J\setminus K \mid \widehat{c}_j=\alpha \text{ and } \widehat{c}_j<0\,\}.
    \end{equation}
    If $K_{\rm gen}=\emptyset$ (i.e., if $\alpha\ge0$ or $K=J$), terminate the iteration.
    
    \State \textbf{Step 2d:} Update the column set by setting
    \begin{equation}
    K\gets K\cup K_{\rm gen}.
    \end{equation} 
\EndWhile

\State \textbf{Step 3:} Finally, solve the following the RMP as a 0-1 integer programming problem using the final column set $K$:
\begin{equation}
\begin{split}
\bm{x}^* =\ & \underset{\bm{x}_K}{\operatorname{argmin}}\; \bm{c}_K^\top\, \bm{x}_K,\\[1ex]
& \text{subject to}\quad A[:,K]\, \bm{x}_K=\bm{1}_M.
\end{split}
\end{equation}

\State \textbf{Return:} The solution $\bm{x}^*$.
\end{algorithmic}
\end{algorithm}

\clearpage

\section{Proposed Method: Column Generation and Elimination Method}\label{sec:ProposedMethod}

In this section, the proposed Column Generation and Elimination method is described. 
This approach improves both the exact solution rate and the accuracy of approximate solutions by relaxing the Column Generation conditions found in conventional methods. At the same time, it introduces a mechanism to eliminate unnecessary columns, which helps to control the growth of the problem size. The pseudocode for obtaining a final solution using this method is presented in {\bf Algorithm}~\ref{alg:ColGenEli}. 
It is assumed that all cost coefficients $c_i$ satisfy $0 < c_i \leq 1$.

The differences between the proposed method and the conventional Column Generation method are summarized as follows.
These four differences were determined based on numerical experiments conducted under various settings: 
\begin{enumerate}
\item In conventional Column Generation, the dual of the linear relaxation of the Set Partitioning Problem is formulated as a linear program with unconstrained continuous variables (see Equation~\eqref{DRMP1}). In contrast, the Column Generation and Elimination Method imposes bounds on these continuous variables, restricting them to the interval $[0,1]$ (see Equation~\eqref{DRMP2}).
\item For the decision of which new columns to generate, the reduced cost (see Equation~\eqref{RC2}) is computed for candidate columns $ J\setminus L $. Instead of selecting only those with negative reduced costs, a relaxed condition that allows non-positive values is adopted (see Equation~\eqref{eq:gencondi}).
\item For the columns already included in the  RMP, the reduced cost (see Equation~\eqref{RC2}) is similarly calculated, and any column with a value of at least 1 is removed (see Equation~\eqref{eq:elicondi}). Removing these columns preserves the optimality of the RMP solution. 
In contrast to conventional Column Generation methods, where columns with negative reduced costs are added to improve the RMP solution, this approach removes redundant columns without sacrificing solution quality. 
\item After the iterative process of Column Generation and Elimination has finished, the indices of the feasible solution obtained initially are added to the index set $ K $ (see {\bf Step 3}). This step guarantees that an integer solution is achieved in the final stage.
\end{enumerate}

\begin{algorithm}
\caption{Column Generation and Elimination Method}
\label{alg:ColGenEli}
\begin{algorithmic}[1]
\State \textbf{Step 1:} Use a greedy algorithm to obtain a feasible solution of the original Set Partitioning Problem, $\bm{z}=(z_1,\ldots,z_N)^\top$. From the index set $J=\{1,\ldots,N\}$, define the initial column set as
$K_{\rm ini}:=\{\,j\in J \mid z_j=1\,\}$.
Initialize $K\gets K_{\rm ini}$ and $L\gets K_{\rm ini}$. Here, $K$ is the set used for column selection in the master problem, and $L$ is the record of columns that have been selected.

\While{True}
    \State \textbf{Step 2a:} Solve the improved dual problem of the linear relaxation of the 
    RMP limited to the index set $K$:
    \begin{equation}
    \begin{split}
    \widehat{\bm{y}} =\ & \underset{\bm{y}=(y_1,\ldots,y_M)^\top}{\operatorname{argmax}}\; \sum_{i=1}^{M} y_i,\\[1mm]
    & \text{subject to}\qquad
    A[:,K]^\top\,\bm{y}\le \bm{c}_K;\qquad
    0\le y_i\le 1\quad (\forall\, i).
    \end{split}
    \label{DRMP2}
    \end{equation}
    
    \State \textbf{Step 2b:} Using the dual solution $\widehat{\bm{y}}$, compute  the reduced cost
    \begin{equation}
    \widehat{c}_j:= c_j-\widehat{\bm{y}}^\top A[:,j]\qquad j\in (J\setminus L)\cup K
    \end{equation}
    for each $j$, and set
    \begin{equation}
    \alpha:= \min\{\,\widehat{c}_j \mid j\in J\setminus L\,\},\quad \beta:= \max\{\,\widehat{c}_j \mid j\in K\,\}.
    \label{RC2}
    \end{equation}

    \State \textbf{Step 2c (Iteration Termination Check):} Define the candidate sets for Column Generation and Elimination as
    \begin{align}
    K_{\rm gen} & := \{\,j\in J\setminus L \mid \widehat{c}_j=\alpha \text{ and } \widehat{c}_j\le 0\,\},
    \label{eq:gencondi}\\
    K_{\rm eli} & := \{\,j\in K \mid \widehat{c}_j=\beta \text{ and } \widehat{c}_j\ge 1\,\}.
    \label{eq:elicondi}
    \end{align}
    If $K_{\rm gen}=\emptyset$ and $K_{\rm eli}=\emptyset$ (i.e., if $\alpha>0$ or $L=J$ and $\beta<1$), exit the iteration.
    
    \State \textbf{Step 2d:} Update the column sets:
    \begin{equation}
    K\gets (K\cup K_{\rm gen})\setminus K_{\rm eli},\quad
    L\gets L\cup K_{\rm gen}.
    \end{equation}    
    
\EndWhile
\State \textbf{Step 3:} Add the initial column set $K_{\rm ini}$ back to $K$; that is, set $K\gets K\cup K_{\rm ini}$.
\State \textbf{Step 4:} Solve the Set Partitioning Problem with variables restricted to the index set $K$ (as in Step 3 of the Column Generation method) to obtain the final solution.
\end{algorithmic}
\end{algorithm}

\clearpage

\section{Evaluation Experiments}\label{sec:Evaluation}

This section describes the evaluation experiments conducted to assess the effectiveness of the proposed method. First, experiments were carried out on a simple Set Partitioning Problem—formulated as a minimum crew scheduling problem with constant cost coefficients (all equal to 1)—using actual bus timetables. Next, experiments were performed on a more complex version of the Set Partitioning Problem, in which the cost coefficients were replaced with random positive numbers less than or equal to  1.

\subsection{Experimental Data Configuration}

In public transportation management, determining the minimum number of crew members required based on established timetables for buses, trains, and similar services is a critical task. This challenge, known as the minimum crew scheduling problem, can be formulated as a Set Partitioning Problem with all cost coefficients equal to 1.

Open data in the GTFS (General Transit Feed Specification) format, representing the bus timetables of a Japanese city, was utilized~\cite{GTFS}. 
From the total of 19 bus routes, subsets consisting of $k \ (2\leq k \leq 11) $ routes were randomly selected ten times to verify the effectiveness of the proposed method at various scales, resulting in a total of $100\ (10\times 10)$ experimental patterns.

\begin{figure}[htbp]
    \centering
    \begin{subfigure}[b]{0.48\linewidth}
        \centering
        \includegraphics[width=\linewidth]{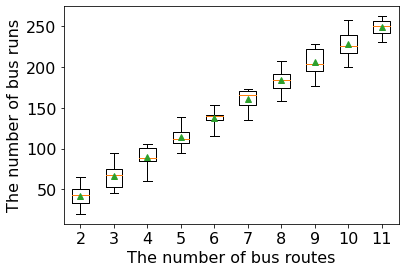}
        \caption{Number of Bus Runs}
    \end{subfigure}
    \hfill
    \begin{subfigure}[b]{0.48\linewidth}
        \centering
        \includegraphics[width=\linewidth]{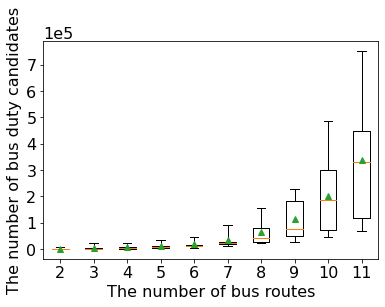}
        \caption{Number of Bus Duty Candidates}
    \end{subfigure}
    \caption{Boxplots showing the number of bus runs (a) and bus duty candidates (b) summarized by route count. In (a), the number of trips increases proportionally with the number of routes. In (b), the number of duty candidates grows exponentially with the route count, accompanied by a larger variation range.}\label{fig:BoxPlot}
\end{figure}

For each selected subset of routes, candidate bus driver duties (complete driving schedules from departure to return at the depot) were generated by exhaustively enumerating all feasible combinations. This enumeration considered constraints such as bus connection requirements and the maximum allowable duty duration (9 hours), and was conducted computationally. The numbers of bus trips and duty candidates for each route are illustrated in a boxplot (Fig.~\ref{fig:BoxPlot}).

\subsection{Optimization Problem Formulation}

Using the notation introduced in Subsection~\ref{subsec:setpartition}, $V$ represents the set of all bus trips, and each $U_i$ corresponds to a candidate duty—namely, a set of bus trips that a single driver can cover within a given duty period.

For each experimental pattern, two types of Set Partitioning Problems were formulated to compare the evaluation results of the conventional method (Column Generation) and the proposed method (Column Generation and Elimination):
\begin{description}
  \item[Problem A] In the case of constant cost coefficients, all cost coefficients are set to $c_i=1$. This formulation corresponds to the minimum crew scheduling problem.
  \item[Problem B] In the case of random cost coefficients, each $c_i$ is randomly selected from the set $\{0.1\cdot k \mid k \text{ is an integer},\, 1\le k\le 10\}$.
\end{description}

\subsection{Evaluation Criteria}

In this paper, the following metrics are used for evaluation.
First, the exact optimal solution for each experimental problem is obtained beforehand using a branch-and-cut algorithm, and its optimal value is defined as the true optimal value $A$.  
\begin{description}
\item[Exact Solutions count (ES count)] 
  The number of times a solution with an objective function value equal to the true optimal value ($A = B$) is obtained. Solutions are considered exact if their objective values match, even if their solution vectors differ.
\item[Approximation Ratio (AR)]
  Defined by the following equation:
  \begin{equation}
  \text{Approximation ratio} := \frac{B - A}{A}
  \end{equation}
A smaller value of AR (closer to zero) indicates higher accuracy of the obtained solution.
\end{description}

\subsection{Experimental Results}\label{subsec:Evaluation}

For 100 experimental patterns, the number of exact solutions obtained and the approximation ratio were evaluated for both Problem A and Problem B. Total and average values, aggregated by the number of routes, are presented in Table~\ref{tbl:result1} and Table~\ref{tbl:result2}.

In Problem A (constant cost coefficients), the conventional method (Column Generation) achieved an exact solution rate of only 3\%, whereas the proposed method (Column Generation and Elimination) improved this rate to 99\%. In contrast, for Problem B (random cost coefficients), both methods tended to show lower exact solution rates and approximation ratios compared to Problem A; however, the proposed method demonstrated a higher improvement in the approximation ratio. Notably, while the approximation ratio of the conventional method remained nearly constant regardless of problem size, the proposed method tended to improve as the problem size increased. 
In particular, a significant improvement in the exact solution rate was observed when the number of bus routes reached its maximum of 11. Similar improvements may be expected for large-scale problems involving several million of variables. 
These results indicate that, in terms of solution accuracy, the proposed method is advantageous when the cost coefficients are uniformly distributed and the problem size increases. However, the precise influence of the cost coefficient distribution and problem size on the evaluation results remains unclear.

\begin{table}[t]
\caption{Evaluation results for Problem A}\label{tbl:result1}
\centering
\begin{tabular*}{12cm}{@{\extracolsep{\fill}}cc|cc|cc}
\multirow{3}{*}{\begin{tabular}{c}\# of  \\routes\end{tabular}}  
& \multirow{3}{*}{\begin{tabular}{c}\# of \\trials\end{tabular}} & 
\multicolumn{2}{c|}{Conventional method} & \multicolumn{2}{c}{\bf Proposed method} \\ 
&  & \multirow{2}{*}{\begin{tabular}{c}ES \\count\end{tabular}} & 
\multirow{2}{*}{\begin{tabular}{c}Average \\of AR \end{tabular}} & 
\multirow{2}{*}{\begin{tabular}{c}ES \\count\end{tabular}} & 
\multirow{2}{*}{\begin{tabular}{c}Average \\of AR \end{tabular}} \\
& & & & & \\ \hline
 2 & 10 & 2 & 0.18 & 10 & 0.00 \\
3 & 10 & 1 & 0.15 & 10 & 0.00 \\
4 & 10 & 0 & 0.22 & 10 & 0.00 \\
5 & 10 & 0 & 0.21 & 10 & 0.00 \\
6 & 10 & 0 & 0.20 & 9 & 0.00 \\
7 & 10 & 0 & 0.18 & 10 & 0.00 \\
8 & 10 & 0 & 0.17 & 10 & 0.00 \\
9 & 10 & 0 & 0.19 & 10 & 0.00 \\
10 & 10 & 0 & 0.16 & 10 & 0.00 \\
11 & 10 & 0 & 0.16 & 10 & 0.00 \\ \hline
Average & 10 & 0.30	& 0.182 & 9.90 & 0.000 \\ 
\end{tabular*}\linebreak\linebreak
\caption{Evaluation results for Problem B}\label{tbl:result2}
\centering
\begin{tabular*}{12cm}{@{\extracolsep{\fill}}cc|cc|cc}
\multirow{3}{*}{\begin{tabular}{c}\# of  \\routes\end{tabular}}  
& \multirow{3}{*}{\begin{tabular}{c}\# of \\trials\end{tabular}} & 
\multicolumn{2}{c|}{Conventional method} & \multicolumn{2}{c}{\bf Proposed method} \\ 
&  & \multirow{2}{*}{\begin{tabular}{c}ES \\count\end{tabular}} & 
\multirow{2}{*}{\begin{tabular}{c}Average \\of AR \end{tabular}} & 
\multirow{2}{*}{\begin{tabular}{c}ES \\count\end{tabular}} & 
\multirow{2}{*}{\begin{tabular}{c}Average \\of AR \end{tabular}} \\
& & & & & \\ \hline
2 & 10 & 2 & 0.33 & 2 & 0.27 \\
3 & 10 & 0 & 0.37 & 0 & 0.29 \\
4 & 10 & 0 & 0.41 & 0 & 0.28 \\
5 & 10 & 0 & 0.47 & 0 & 0.28 \\
6 & 10 & 0 & 0.41 & 0 & 0.20 \\
7 & 10 & 0 & 0.39 & 0 & 0.19 \\
8 & 10 & 0 & 0.39 & 0 & 0.11 \\
9 & 10 & 0 & 0.42 & 0 & 0.10 \\
10 & 10 & 0 & 0.41 & 1 & 0.04 \\
11 & 10 & 0 & 0.39 & 4 & 0.02 \\ \hline
Average & 10 & 0.20 & 0.399 & 0.70	& 0.178 \\ 
\end{tabular*}
\end{table}

In addition, to monitor the status of the iterative calculation process, the average number of iterations and the compression ratio, by route count, were computed for both Problem A and Problem B. The results are shown in Table~\ref{tbl:result3} and Table~\ref{tbl:result4}. The following definitions apply:

\begin{description}
\item[Iteration Count] Refers to the number of iterations required for convergence of Steps 2a–2d in Algorithm~\ref{alg:ColGen} and Algorithm~\ref{alg:ColGenEli}.
\item[Compression Ratio] Defined as
\begin{equation}
\text{Compression Ratio} := \frac{\#K}{\#J},
\end{equation}
where $\#X$ denotes the number of elements in set $X$.
\end{description}

In Problem A, the proposed method exhibited higher iteration counts and compression ratios than the conventional method. Although the proposed method relaxes the Column Generation conditions compared to the conventional method, the relatively small number of eliminated columns necessitates more iterations for convergence, resulting in a higher compression ratio (i.e., a larger value). Moreover, while the compression ratio of the conventional method tended to decrease with increasing problem size, that of the proposed method remained nearly constant. In contrast, for Problem B the difference between the two methods diminished. Furthermore, as the problem size increased, the proposed method also tended to exhibit a lower compression ratio (i.e., a higher degree of compression). These results suggest that in terms of computational cost, the difference between the conventional and proposed methods narrows when the cost coefficients exhibit high variance and the problem size increases. However, the reason why the variance of cost coefficients affects the evaluation results remains unclear.

\begin{table}[t]
\caption{Status of the computational process in Problem A}\label{tbl:result3}
\centering
\begin{tabular*}{12cm}{@{\extracolsep{\fill}}c|cc|cc}
\multirow{3}{*}{\begin{tabular}{c}\# of  \\routes\end{tabular}}   & 
\multicolumn{2}{c|}{Conventional method} & \multicolumn{2}{c}{\bf Proposed method} \\ 
&  \multirow{2}{*}{\begin{tabular}{c}Iteration \\count\end{tabular}} & 
\multirow{2}{*}{\begin{tabular}{c}Compression \\ratio \end{tabular}} & 
\multirow{2}{*}{\begin{tabular}{c}Iteration \\count\end{tabular}} & 
\multirow{2}{*}{\begin{tabular}{c}Compression \\ratio \end{tabular}} \\ 
& & & &  \\ \hline
2 & 4.6 & 0.15 & 11.5 & 0.36 \\
3 & 5.8 & 0.13 & 12.4 & 0.37 \\
4 & 6.8 & 0.13 & 16.0 & 0.33 \\
5 & 9.6 & 0.10 & 18.1 & 0.32 \\
6 & 11.5 & 0.07 & 22.2 & 0.34 \\
7 & 12.7 & 0.07 & 23.4 & 0.37 \\
8 & 12.3 & 0.05 & 24.6 & 0.36 \\
9 & 13.8 & 0.05 & 32.5 & 0.36 \\
10 & 16.2 & 0.06 & 33.1 & 0.35 \\
11 & 16.8 & 0.04 & 43.3 & 0.37 \\ \hline
Average  & 11.0 & 0.085 & 23.7 & 0.353 \\ 
\end{tabular*}\linebreak\linebreak
\caption{Status of the computational process in Problem B}\label{tbl:result4}
\centering
\begin{tabular*}{12cm}{@{\extracolsep{\fill}}c|cc|cc}
\multirow{3}{*}{\begin{tabular}{c}\# of  \\routes\end{tabular}}   & 
\multicolumn{2}{c|}{Conventional method} & \multicolumn{2}{c}{\bf Proposed method} \\ 
&  \multirow{2}{*}{\begin{tabular}{c}Iteration \\count\end{tabular}} & 
\multirow{2}{*}{\begin{tabular}{c}Compression \\ratio \end{tabular}} & 
\multirow{2}{*}{\begin{tabular}{c}Iteration \\count\end{tabular}} & 
\multirow{2}{*}{\begin{tabular}{c}Compression \\ratio \end{tabular}} \\ 
& & & &  \\ \hline
2 & 11.8 & 0.08 & 16.5 & 0.11 \\
3 & 13.8 & 0.05 & 15.7 & 0.08 \\
4 & 15.6 & 0.04 & 22.2 & 0.07 \\
5 & 22.0 & 0.03 & 22.8 & 0.06 \\
6 & 26.8 & 0.02 & 29.8 & 0.06 \\
7 & 25.0 & 0.02 & 30.1 & 0.06 \\
8 & 30.6 & 0.02 & 33.8 & 0.05 \\
9 & 29.7 & 0.01 & 32.2 & 0.05 \\
10 & 34.2 & 0.01 & 37.9 & 0.06 \\
11 & 36.4 & 0.01 & 40.3 & 0.06 \\ \hline
Average  & 24.6 & 0.029 & 28.1 & 0.066 \\ 
\end{tabular*}
\end{table}

\clearpage

\section{Conclusions}\label{sec:Conc}

A novel approach, termed the Column Generation and Elimination Method, has been introduced for solving large-scale Set Partitioning Problems. This method improves both the exact solution rate and the accuracy of approximate solutions by relaxing the Column Generation conditions, while significantly reducing the problem size through a mechanism for eliminating redundant columns. 
Moreover, since the proposed method does not require branching operations, its implementation remains as straightforward as the conventional Column Generation method. 
Numerical experiments demonstrated the superiority of the proposed method in terms of solution quality using open GTFS data of real-world bus routes. 
In particular, when cost coefficients are nearly uniform and the problem size increases, the proposed method tends to yield more accurate solutions. Conversely, regarding computational cost, the performance gap between conventional and proposed methods diminishes when the variance of cost coefficients is high and the problem scale expands. 
However, a theoretical analysis explaining the advantages of the proposed method remains an open challenge for future studies.

Existing studies on Column Generation include efforts to accelerate the process using machine learning techniques~\cite{MachineLearning1,MachineLearning2,MachineLearning3} and investigations into hybrid approaches that combine it with other methods for faster computation~\cite{CombiningOtherwise1}. Future research will explore further acceleration and improvements in solution accuracy by building upon these approaches.

\clearpage

\appendix

\section{Greedy Algorithm}\label{sec:greedy}

\begin{algorithm}[H]
\caption{Greedy Algorithm}
\label{alg:greedy}
\begin{algorithmic}[1]
\State \textbf{Step 1:} Given a set $V = \{v_i\}_{i\in I}$ and a family of its subsets $\{U_j\}_{j\in J}$, where $I$ and $J$ are index sets. For each $j\in J$, let $c_j > 0$ denote the cost of selecting subset $U_j$ to cover $V$. Define $K$ as the index set for the subsets chosen to cover $V$, and initialize $K := \emptyset$.

\While{True}
    \State \textbf{Step 2a:} Define the index set $L$ of subsets that cover portions of $V$ not yet covered by the selected subsets:
    \begin{equation}
    L := \{\,j\in J \mid U_j \cap U_k = \emptyset \quad \forall\, k\in K\,\}.
    \end{equation}
    \State \textbf{Step 2b (Iteration Termination Check):} If $L$ is empty, terminate the iteration.
    
    \State \textbf{Step 2c:} For each $j\in L$, denote the number of elements in $U_j$ by $\#U_j$. Then, select
    \begin{equation}
    \widehat{j} := \underset{j\in L}{\operatorname{argmax}}\ \frac{\#U_j}{c_j}.
    \end{equation}
    \State \textbf{Step 2d:} Add the selected index $\widehat{j}$ to $K$:
    \begin{equation}
    K \gets K \cup \{\widehat{j}\}.
    \end{equation}
\EndWhile
\State \textbf{Step 3:} For each $j\in J$, define
\begin{equation}
x_j =
\begin{cases}
1, & \text{if } j\in K,\\[1ex]
0, & \text{otherwise,}
\end{cases}
\end{equation}
and construct the initial feasible solution as the binary vector $\bm{x}_{\rm ini} := (x_j)_{j\in J}$.
\State \textbf{Return:} $\bm{x}_{\rm ini}$.
\end{algorithmic}
\end{algorithm}


\end{document}